\documentclass[12pt,a4paper]{article}

\usepackage{amssymb}

\newcommand{\stl}{\vspace{3mm}}

\newtheorem{defi}{\sc Definition}[section]
\newtheorem{theo}[defi]{\sc Theorem}
\newtheorem{prop}[defi]{\sc Proposition}
\newtheorem{nota}[defi]{\sc Notation}

\newtheorem{rema}[defi]{\sc Remark}
\newtheorem{lemm}[defi]{\sc Lemma}
\newtheorem{exem}[defi]{\sc Example}

\newenvironment{demo}[1][]{\noindent {\it Proof.} \protect\nopagebreak
       \rm #1}{\protect\nopagebreak $\square $ \par\stl}

\begin{document}

\begin{center}{\bf Intersection homology ${\cal D}$-Module 
and Bernstein~polynomials 

associated with a complete intersection}
\end{center}
\medskip
\begin{center}
{\large Tristan Torrelli}\footnote{Laboratoire 
J.A. Dieudonn\'e, UMR
du CNRS 6621, Universit\'e de Nice Sophia-Antipolis,
Parc Valrose, 06108 Nice Cedex 2, France. {\it E-mail:} 
tristan$\_$torrelli@yahoo.fr

\noindent {\it 2000 Mathematics Subject Classification:} 32S40, 
32C38, 32C40, 32C25, 14B05.

\noindent {\it Keywords:}  intersection holomology 
${\cal D}$-modules, local algebraic cohomology group, 
complete intersections, Bernstein-Sato functional equations.}
\end{center}

\bigskip

\noindent{\sc Abstract}. 
Let $X$ be a complex analytic manifold. Given a closed subspace $Y\subset X$ of pure codimension $p\geq 1$, we consider the sheaf of local algebraic 
cohomology $H^p_{[Y]}({\cal O}_X)$, and ${\cal L}(Y,X)\subset H^p_{[Y]}({\cal O}_X)$ 
 the intersection homology ${\cal D}_X$-Module of Brylinski-Kashiwara.
 We give here an algebraic characterization of the spaces $Y$ such that
${\cal L}(Y,X)$ coincides with $H^p_{[Y]}({\cal O}_X)$, in terms of Bernstein-Sato functional equations.

\medskip

\section{Introduction}

Let $X$ be a complex analytic manifold of dimension $n\geq 2$, 
${\cal O}_X$ be the sheaf of holomorphic functions on $X$ and
${\cal D}_X$ the sheaf of differential operators with holomorphic
coefficients. At a point $x\in X$, we identify the stalk 
${\cal O}_{X,x}$ (resp. ${\cal D}_{X,x}$) with the ring 
${\cal O}={\bf C}\{x_1,\ldots,x_n\}$ (resp. 
${\cal D}={\cal O}\langle\partial/\partial x_1,\ldots,
\partial/\partial x_n\rangle$).

Given a closed subspace $Y\subset X$ of pure codimension $p\geq 1$, 
we denote by $H^p_{[Y]}({\cal O}_X)$ the sheaf of local algebraic 
cohomology with support in $Y$. Let 
${\cal L}(Y,X)\subset H^p_{[Y]}({\cal O}_X)$ be the intersection
homology ${\cal D}_X$-Module of Brylinski-Kashiwara (\cite{Br}). This is
the smallest ${\cal D}_X$-submodule of $H^p_{[Y]}({\cal O}_X)$
which coincides with $H^p_{[Y]}({\cal O}_X)$ at the generic points
of $Y$ (\cite{Br}, \cite{BK}).

\stl

A natural problem is to characterize the subspaces $Y$ such that
${\cal L}(Y,X)$ coincides with $H^p_{[Y]}({\cal O}_X)$. We prove here 
that it may be done locally using Bernstein functional equations. This
supplements a work of D. Massey (\cite{Mas}), who studies
the analogous problem with a topological viewpoint. Indeed,
from the Riemann-Hilbert 
correspondence of Kashiwara-Mebkhout (\cite{K4}, \cite{Mb2}),
the regular holonomic ${\cal D}_X$-Module $H^p_{[Y]}({\cal O}_X)$ 
corresponds to the perverse sheaf ${\bf C}_Y[n-p]$ (\cite{G}, \cite{K2}, 
\cite{Mb1}) where as ${\cal L}(Y,X)$ corresponds to the
intersection complex $IC^{\bullet}_Y$ (\cite{Br}). By this way, this
condition ${\cal L}(Y,X)=H^p_{[Y]}({\cal O}_X)$ is equivalent to
the following one: \textsl{the real link of $Y$ at a point $x\in Y$ is
a rational homology sphere}.

\stl

 First of all, we
have  an explicit local description of ${\cal L}(Y,X)$. This comes
from the following result, due to D. Barlet and M. Kashiwara. 

\begin{theo}[\cite{BK}]
  The fundamental class 
$C^Y_X\in H^p_{[Y]}({\cal O}_X)\otimes \Omega^p_X$ of $Y$ in $X$ 
belongs to ${\cal L}(Y,X)\otimes \Omega^p_X$.
\end{theo}  

For more details about $C^Y_X$, see \cite{B1}. In particular, if $h$   
is an analytic morphism 
$(h_1,\ldots, h_p):(X,x)\rightarrow({\bf C}^p,0)$ 
 which defines the complete 
intersection $(Y,x)$ - reduced or not -, then the inclusion 
${\cal L}(Y,X)_x\subset H^p_{[Y]}({\cal O}_X)_x$ may be identified with:
$${\cal L}_h=\sum_{1\leq k_1<\cdots < k_p\leq n}{\cal D}\cdot
\frac{\stackrel{.}{m_{k_1,\ldots,k_p}(h)}}{h_1\cdots h_p} \subset 
{\cal R}_h=\frac{{\cal O}[1/h_1\cdots h_p]}
{\sum_{i=1}^p{\cal O}[1/h_1\cdots \check{h}_i\cdots h_p]}$$
where $m_{k_1,\ldots,k_p}(h)\in{\cal O}$ is the determinant
of the columns $k_1,\ldots, k_p$ of the Jacobian matrix of $h$. 
In the following, ${\cal J}_h\subset{\cal O}$ denotes the ideal 
generated by the $m_{k_1,\ldots,k_p}(h)$, and
$\delta_h\in{\cal R}_h$ the section defined by $1/h_1\cdots h_p$.

\stl

When $Y$ is a hypersurface, we have the following characterization. 

\begin{theo} \label{LaiR}
 Let $Y\subset X$ be a hypersurface and $h\in{\cal O}_{X,x}$ denote 
a local equation of $Y$ at a point $x\in Y$. 
The following conditions are equivalent:
 \begin{enumerate}  
   \item ${\mathcal L}(Y,X)_x$ coincides with
   $H^p_{[Y]}({\cal O}_X)_x$.
   \item The reduced Bernstein
    polynomial of $h$ has no integral root.
   \item $1$ is not an eigenvalue of the monodromy acting on the
    reduced cohomology of the fibers of the Milnor fibrations of $h$
    around any singular points of $Y$ contained in some open neighborhood of $x$ in $Y$. 
 \end{enumerate}
\end{theo}

Let us recall that the {\em Bernstein polynomial} $b_f(s)$ of a 
nonzero germ $f\in{\cal O}$ is the monic generator of the ideal
of the polynomials $b(s)\in{\bf C}[s]$ such that:
$$b(s)f^s=P(s)\cdot f^{s+1}$$
in ${\cal O}[1/f,s]f^s$, where $P(s)\in{\cal D}[s]={\cal D}
\otimes {\bf C}[s]$. The existence of such a nontrivial
equation was proved by M. Kashiwara (\cite{11}).
 When $f$ is not a unit, it is easy to check that 
$-1$ is a root of $b_f(s)$. The quotient of $b_f(s)$ by
$(s+1)$ is the so-called {\em reduced Bernstein polynomial}
of $f$, denoted $\tilde{b}_f(s)$. Let us recall that their roots are
rational negative numbers in $]-n,0[$ (see \cite{Sait} for the general case, 
\cite{Var} for the isolated singularity case).

\begin{exem}{\em \label{exbern}
  Let $f=x_1^2+\cdots+x_n^2$. It is easy to prove that $b_f(s)$
is equal to $(s+1)(s+n/2)$, by using the identity:
$$\left[\frac{\partial}{\partial x_1} ^2+\cdots +\frac{\partial}{\partial x_n} ^2\right]\cdot
f^{s+1}=(2s+2)(2s+n)f^s$$ In particular, ${\cal R}_f$ coincides with
${\cal L}_f$ if and only if $n$ is odd.} 
\end{exem}

These polynomials are famous because of the link of their
roots with the monodromy of the Milnor fibration associated
with $f$. This was established by B. Malgrange \cite{M2} and M. 
Kashiwara \cite{K3}. More generally, by using the algebraic microlocalization, 
M. Saito \cite{Sait} prove that 
$\{ e^{-2i\pi\alpha}\,|\, \alpha\mbox{ root of }\tilde{b}_f(s)\}$ 
coincides with the set of the eigenvalues of the monodromy acting
on the Grothendieck-Deligne vanishing cycle sheaf $\phi_f{\bf C}_{X(x)}$ 
(where $X(x)\subset X$ is a sufficiently small neighborhood of $x$).
Thus the equivalence $2\Leftrightarrow 3$ is an easy consequence of 
this deep fact.

\stl

We give a direct proof of $1\Leftrightarrow 2$ in part 4.

\begin{rema}{\em
 In \cite{B2}, D. Barlet gives a characterization of 
$3$ in terms of the meromorphic continuation of the current 
$\int_{X(x)} f^\lambda\square$.}
\end{rema}

\begin{rema}{\em
The equivalence $1\Leftrightarrow 3$ for the isolated singularity case may be due
to J. Milnor \cite{Mil} using the Wang sequence relating the cohomology of the link with the 
Milnor cohomology. In general, this equivalence is well-known to specialists. It can be proved by using a formalism of
weights and by reducing it to the assertion that the $N$-primitive part of the middle graded
piece of the monodromy weight filtration on the nearby cycle sheaf is the intersection complex
(this last assertion is proved in \cite{Sait90} (4.5.8) for instance). It would be quite interesting if one can prove the equivalence between 1 and 3 by using only the theory of $\mathcal D$-modules.
}
\end{rema}

In the case of hypersurfaces, it is well known that 
condition $1$ requires a 
strong kind of irreducibility. This may be refinded in terms of 
Bernstein polynomial. 

\begin{prop}\label{redBred}
   Let $f\in{\cal O}$ be a nonzero germ such that $f(0)=0$. Assume that the origin belongs
   to the closure of the points where $f$ is locally reducible. 
   Then $-1$ is a root of the reduced Bernstein polynomial of $f$.
\end{prop}

\begin{exem}{\em \label{ExLnotR}
If $f=x_1^2+x_3x_2^2$, then $(s+1)^2$
divides $b_f(s)$ because $f^{-1}\{0\}\subset$~${\bf C}^3$ is reducible at
any $(0,0,\lambda)$, $\lambda\not=0$ (in fact, we have: 
$b_f(s)=(s+1)^2(s+3/2)$).}
\end{exem}

What may be done in higher codimensions ? If 
$f\in{\cal O}$ is such that $(h,f)$ defines a complete intersection,
we can consider the Bernstein polynomial $b_f(\delta_h,s)$ of $f$
associated with $\delta_h\in{\cal R}_h$. Indeed, we again have 
nontrivial functional equations:
$$b(s)\delta_h f^s=P(s)\cdot \delta_h f^{s+1}$$
with $P(s)\in{\cal D}[s]$ (see part \ref{PartBern}). This 
polynomial $b_f(\delta_h,s)$ 
is again a multiple of $(s+1)$, and we can define a reduced Bernstein
polynomial $\tilde{b}_f(\delta_h,s)$ as above. Meanwhile, in order 
to generalize Theorem \ref{LaiR}, the suitable Bernstein polynomial
 is neither $\tilde{b}_f(\delta_h,s)$ nor $b_f(\delta_h,s)$, but a third
one trapped between these two.

\begin{nota} 
Given a morphism 
$(h,f)=(h_1,\ldots,h_p,f):({\bf C}^n,0)\rightarrow ({\bf C}^{p+1},0)$
defining a complete intersection, we denote by
$b'_f(h,s)$ the monic generator of the ideal of
polynomials $b(s)\in{\bf C}[s]$ such that:
\begin{eqnarray} \label{eqcarbprim}
 b(s)\delta_h f^s &\in & {\cal D}[s]({\cal J}_{h,f},f)\delta_hf^s\ .
\end{eqnarray}
\end{nota}

\begin{lemm} \label{relatc}
 The polynomial $b'_f(h,s)$ divides $b_f(\delta_h,s)$,
and $b_f(\delta_h,s)$ divides $(s+1)b'_f(h,s)$. In other words, 
$b'_f(h,s)$ is either $b_f(\delta_h,s)$ or $b_f(\delta_h,s)/(s+1)$.
\end{lemm}

The first assertion is clear since 
${\cal D}[s]\delta_h f^{s+1}\subset {\cal D}[s]({\cal J}_{h,f},f)\delta_h f^s$.
The second relation uses the identities: 
\begin{eqnarray*}
(s+1)m_{k_1,\ldots,k_{p+1}}(h,f)\delta_hf^s & = &
\underbrace{\left[ \sum_{i=1}^{p+1}(-1)^{p+i+1} 
m_{k_1,\ldots,\check{k}_i,\ldots,k_{p+1}}(h)
\frac{\partial}{\partial x_{k_i}}\right]}_{\Delta^h_{k_1,\ldots, k_{p+1}}} 
 \cdot\, \delta_hf^{s+1}
\end{eqnarray*}
for $1\leq k_1<\cdots <k_{p+1}\leq n$, where the vector field 
$\Delta^h_{k_1,\ldots, k_{p+1}}$ annihilates $\delta_h$. 
In particular, we have: $(s+1)b'_f(h,s)\delta_hf^s\in {\cal D}[s]\delta_h f^{s+1}$,
and the assertion follows.

\stl

As a consequence of this result, $b'_f(h,s)$ coincides with $\tilde{b}_f(\delta_h,s)$
when $-1$ is not a root of $b'_f(h,s)$; but it is not always true
(see part \ref{partcompa}). We point out some facts about this polynomial
in part \ref{partcompa}.

\begin{theo} \label{LaiRgen}
  Let $Y\subset X$ be a closed subspace of pure codimension
$p+1\geq 2$, and $x\in Y$. Let
$(h,f)=(h_1,\ldots,h_p,f):(X,x)\rightarrow ({\bf C}^{p+1},0)$
be an analytic morphism such that the common zero set of 
$h_1,\ldots,h_p,f$ is $Y$ in a neighbourhood of $x$. Up to replace $h_i$
by $h_i^m$ for some non negative integer $m\geq 1$, let us assume that
${\cal D}\delta_{h}={\cal R}_h$. The
following conditions are equivalent:
 \begin{enumerate}  
   \item ${\mathcal L}(Y,X)_x$ coincides with $H^p_{[Y]}({\cal O}_X)_x$.
   \item The polynomial $b'_f(h,s)$ has no strictly 
    negative integral root.
 \end{enumerate}
\end{theo}

\stl

Let us observe that the condition ${\cal D}\delta_{h}={\cal R}_h$
is not at all a constraining condition on $(Y,x)$. Moreover, using the
boundaries of the roots of the classical Bernstein polynomial, on can
take $m=n-1$ (since $1/(h_1\cdots h_p)^{n-1}$ generates the 
$\mathcal D$-module ${\mathcal O}[1/h_1\cdots h_p]$, using Proposition 
\ref{propmhs} below). Finally, one can observe that this 
technical condition ${\cal D}\delta_{h}={\cal R}_h$ is difficult
to verify in practice. Thus, let us give an inductive criterion.  

\begin{prop} \label{criterion}
 Let $h=(h_1,\ldots,h_p):(X,x)\rightarrow ({\bf C}^p,0)$
be an analytic morphism defining a germ of complete
intersection of codimension $p\geq 1$. Assume that $-1$ is the only integral 
root of the Bernstein polynomial $b_{h_1}(s)$.
 Moreover, if $p\geq 2$, assume that $-1$ is the smallest integral
 root of $b_{h_{i+1}}(\delta_{\tilde{h}_i},s)$ 
 with $\tilde{h}_i=(h_1,\ldots,h_i):(X,x)\rightarrow ({\bf C}^i,0)$, for $1\leq i\leq p-1$. 
 Then the left $\mathcal D$-module ${\cal R}_h$ is generated by $\delta_h$.   
\end{prop}

\begin{exem} {\em
Let $n=3$, $p=2$, $h_1=x_1^2+x_2^3+x_3^4$ and $h_2=x_1^2-x_2^3+2x_3^4$. As 
$h_1$ defines an isolated singularity and $h=(h_1,h_2)$ defines a weighted-homogeneous 
complete intersection isolated singularity, we have closed formulas for
$b_{h_1}(s)$ and $b_{h_2}(\delta_{h_1},s)$, see \cite{Ya}, \cite{Cras1}. From 
the explicit expression of these two polynomials, we see that they have no integral 
root smaller than $-1$. Thus $\delta_h$ generates ${\mathcal R}_h$.}
\end{exem}

The proofs of Theorems \ref{LaiR} \& \ref{LaiRgen} are given in 
part \ref{proof}. They are based
on a natural generalization of a classical result due to M.
Kashiwara which links the roots of $b_f(s)$ to some generators
of ${\cal O}[1/f]f^\alpha$, $\alpha\in{\bf C}$ (Proposition
\ref{propmhs}). The last part is devoted to remarks and comments
about Theorem $\ref{LaiRgen}$.

\stl

\noindent{\bf Aknowledgements.} This research has 
been supported by a Marie Curie Fellowship of the European 
Community (programme FP5, contract HPMD-CT-2001-00097). 
The author is very grateful to the Departamento de \'Algebra, 
Geometr\'{\i}a y Topolog\'{\i}a (Universidad de Valladolid)
 for hospitality during the fellowship, to the Lehrstuhl VI 
f\"ur Mathematik (Universit\"at Mannheim) for hospitality in 
November 2005, and to Morihiko Saito for judicious comments.

\section{Bernstein polynomials associated with a section of a holonomic
$\cal D$-module}
\label{PartBern}
In this paragraph, we recall some results about Bernstein polynomials
associated with a section of a holonomic ${\cal D}_X$-Module.

\stl

Given a nonzero germ $f\in {\mathcal O}_{X,x}\cong {\cal O}$ 
and a local section $m\in {\mathcal M}_x$ of a holonomic
${\mathcal D}_X$-Module ${\mathcal M}$ without $f$-torsion, 
M. Kashiwara \cite{K2} proved that there exists
a functional equation:
$$b(s)mf^s=P(s)\cdot mf^{s+1}$$
in $({\cal D}m)\otimes{\cal O}[1/f,s]f^s$, where $P(s)\in{\cal D}[s]={\cal D}
\otimes {\bf C}[s]$ and $b(s)\in{\bf C}[s]$ are nonzero. The
{\em Bernstein polynomial} of $f$ associated with $m$, denoted by $b_f(m,s)$,
is the monic generator of the ideal of polynomials 
$b(s)\in{\bf C}[s]$ which satisfies such an equation. 
When $f$ is not a unit, it is easy to check that if
$m\in {\mathcal M}_x- f{\mathcal M}_x$,
 then $-1$ is a root of $b_f(m,s)$.

Of course, if ${\cal M}={\cal O}_X$ and $m=1$, this is the classical
notion recalled in the introduction. 

\stl

Let us recall that when ${\cal M}$ is a regular holonomic 
${\cal D}_X$-Module, the roots of the polynomials $b_f(m,s)$ 
are closely linked to the eigenvalues of the monodromy of the
perverse sheaf $\psi_f(Sol({\cal M}))$ around $x$, the 
Grothendieck-Deligne nearby cycle sheaf, see \cite{Mb3} for example. Here
$Sol({\cal M})$ denotes the complex $RHom_{{\cal D}_X}({\cal M},{\cal O}_X)$
of holomorphic solutions of ${\cal M}$, and the relation is similar to the 
one given in the introduction (since $Sol({\cal O}_X)\cong {\bf C}_X$).
This comes from the algebraic construction of vanishing cycles,
using Malgrange-Kashiwara $V$-filtration \cite{M2}, \cite{K3}.

\stl
  
  Now, if $Y\subset X$ is a subspace of pure codimension $p$,
then the regular holonomic ${\cal D}_X$-Module  $H^p_{[Y]}({\cal O}_X)$
corresponds to ${\bf C}_Y\,[n-p]$. Thus the roots of the polynomials
$b_f(\delta,s)$, $\delta \in H^p_{[Y]}({\cal O}_X)_x$, are 
linked to the monodromy associated with 
$f:(Y,0)\rightarrow({\bf C},0)$. For more results about
these polynomials, see \cite{16}.

\section{The polynomials $b'_f(h,s)$ and $\tilde{b}_f(\delta_h,s)$} \label{partcompa}

Let us recall that $b'_f(h,s)$ is always equal to one
of the two polynomials $b_f(\delta_h,s)$ and $\tilde{b}_f(\delta_h,s)$. 
In this paragraph, we point out some facts about  
 these Bernstein polynomials associated with 
an analytic morphism 
$(h,f)=(h_1,\ldots,h_p,f):({\bf C}^n,0)\rightarrow({\bf C}^{p+1},0)$
defining a complete intersection.

\stl

First we have a closed formula for $b'_f(h,s)$
  when $h$ and $(h,f)$ define weighted-homogeneous isolated
complete intersection singularities.

\begin{prop}[\cite{Cras1}] \label{propc}
   Let $f,h_1,\ldots,h_p\in{\bf C}[x_1,\ldots,x_n]$, $p<n$, be
some weighted-homogeneous of degree $1$, $\rho_1,\ldots,\rho_p
\in{\bf Q}^{*+}$ for a system of weights 
$\alpha=(\alpha_1,\ldots,\alpha_n)\in({\bf Q}^{*+})^n$.
Assume that the morphisms $h=(h_1,\ldots, h_p)$ and $(h,f)$ define
two germs of isolated complete intersection singularities. Then the
polynomial $b'_f(h,s)$ is equal to:
$$\prod_{q\in\Pi}(s+|\alpha|-\rho_h+q)$$
where $|\alpha|=\sum_{i=1}^n\alpha_i$, $\rho_h=\sum_{i=1}^p\rho_i$ and
$\Pi\subset {\bf Q}^{+}$ is the set of the weights of the
elements of a weighted-homogeneous basis of 
${\cal O}/(f,h_1,\ldots, h_p) {\cal O}+{\cal J}_{h,f}$.
\end{prop}

When $h$ is not reduced, the determination of $b'_f(h,s)$ is
more difficult, even if $(h,f)$ is a homogeneous morphism.

\begin{exem}{\em     \label{excalcbprime}
Let $p=1$, $f=x_1$ and $h=(x_1^2+\cdots +x_n^2)^\ell$ with $\ell\geq 1$.
By using a formula given in \cite{18}, Remark 4.12, the polynomial 
 $\tilde{b}_{x_1}(\delta_h,s)$ is equal to $(s+n-2\ell)$ for any 
 $\ell\in{\bf N}^*={\bf C}^* \cap {\bf N}$. For $\ell\geq n/2$, let us determine 
 $b'_{x_1}(h,s)$ with the help of Theorem \ref{LaiRgen}.
From Example \ref{exbern}, we have ${\cal R}_h={\cal D}\delta_h$ if 
$\ell\geq n/2$, and ${\cal L}_{h,x_1}={\cal R}_{h,x_1}$ if and only if
$n$ is even.

 If $n\leq 2\ell$ is odd, $b'_{x_1}(h,s)$ must coincide with
 $b_{x_1}(\delta_h,s)=(s+1)(s+n-2\ell)$ 
because of Theorem \ref{LaiRgen} 
(since ${\cal L}_{h,x_1}\not={\cal R}_{h,x_1}$). 
On the other hand, if $n\leq 2\ell$ is even, we have 
$b'_{x_1}(h,s)=(s+n-2\ell)=\tilde{b}_{x_1}(\delta_{h},s)$ 
 by the same arguments.
Let us refind this last fact by a direct calculus.  

As $\tilde{b}_{x_1}(\delta_h,s)=(s+n-2\ell)$ divides $b'_{x_1}(h,s)$, we just have
to check that this polynomial $(s+n-2\ell)$ provides a functional equation
for $\tilde{b}_{x_1}(\delta_h,s)$ when $n$ is even. First, we observe that
$$(s+n-2\ell)\delta_h x_1^s=\left[\sum_{i=1}^n \frac{\partial}{\partial x_i}x_i \right]\cdot \delta_h x_1^s$$
If $\ell=1$, we get the result (since ${\cal J}_{h,x_1}=(x_2,\dots,x_n){\cal O}$ 
in that case). Now we
assume that $\ell\geq 2$. Let us prove that $x_i\delta_h x_1^s$ belongs to
${\cal D}({\cal J}_{h,x_1},x_1)\delta_h x_1^s$ for $2\leq i\leq n$. We denote
by $g$ the polynomial $x_1^2+\cdots+x_n^2$ and, for 
$0\leq j \leq \ell-1$, by 
${\cal N}_j\subset{\cal D}({\cal J}_{h,x_1},x_1)\delta_h x_1^s$ the submodule
generated by $x_1\delta_h x_1^s$, $x_2g^j\delta_hx_1^s$, $\ldots,$
$x_ng^j\delta_hx_1^s$. In particular, $h=g^\ell$, 
${\cal N}_{\ell-1}={\cal D}({\cal J}_{h,x_1},x_1)\delta_h x_1^s$ and 
${\cal N}_{j+1}\subset{\cal N}_j$ for $1\leq j\leq \ell-2$. To conclude, we have to check
that ${\cal N}_0={\cal N}_{\ell-1}$.

By a direct computation, we obtain the identity:
$$\frac{\partial}{\partial x_i}\left[2(j-\ell)g^{j-1}x_1^2 +\sum_{k=2}^n \frac{\partial}{\partial x_k}x_kg^j \right]\cdot\, \delta_hx_1^s\ =
2(j-\ell)(n+2(j-\ell)-1)x_i g^{j-1}\delta_hx_1^s$$
for $2\leq i\leq n$, $j> 0$. As $n$ is even, we deduce that 
$x_ig^{j-1}\delta_hx_1^s$ belongs to ${\cal N}_j$ for $2\leq i\leq n$.
In other words, ${\cal N}_{j-1}={\cal N}_j$ for $1\leq j\leq\ell-1$; thus
${\cal N}_{0}={\cal N}_{\ell-1}$, as it was expected.
}
\end{exem}

As the polynomial $b'_f(h,s)$ plays the rule of $\tilde{b}_f(s)$ in
Theorem \ref{LaiRgen}, a natural question is to compare these polynomials 
$b'_f(h,s)$ and $\tilde{b}_f(\delta_h,s)$. Of course, when $(s+1)$ is not a
factor of $b'_f(h,s)$, 
then $b'_f(h,s)$ must coincide with $\tilde{b}_f(\delta_h,s)$;
from Theorem \ref{LaiRgen}, this sufficient condition is satisfied when 
${\cal D}\delta_h={\cal R}_h$ and ${\cal R}_{h,f}={\cal L}_{h,f}$.
But in general, all the cases are possible (see Example \ref{excalcbprime}); 
nevertheless, we do not have found an example with $f$ and $h$ reduced and
$b'_f(h,s)=b_f(\delta_h,s)$. Is $b'_f(h,s)$ always equal to $\tilde{b}_f(h,s)$ 
in this context ? The question is open.
Let us study this problem when $(h,f)$ defines an isolated complete
intersection singularity. In that case, let us consider the short exact
sequence:

$$0\rightarrow {\cal K}\hookrightarrow \frac{{\cal D}[s]\delta_h f^s}
{{\cal D}[s]({\cal J}_{h,f},f)\delta_h f^s} \twoheadrightarrow 
(s+1)\frac{{\cal D}[s]\delta_h f^s}{{\cal D}[s]\delta_h f^{s+1}}\rightarrow 0 $$
where the three $\cal D$-modules are supported by the origin.
Thus the polynomial $b'_f(h,s)$ is equal to 
$\mbox{l.c.m}(s+1,\tilde{b}_f(\delta_h,s))$ if ${\cal K}\not=0$
and it coincides with $\tilde{b}_f(\delta_h,s)$ if not. Remark that
$\cal K$ is not very explicit, since there does not exist a general 
Bernstein functional equation which defines $\tilde{b}_f(\delta_h,s)$ - 
contrarily to $\tilde{b}_f(s)$, see part \ref{proof}. In \cite{T1}, \cite{18},
we have investigated some contexts where such a functional equation
may be given. In particular, this may be done when the following 
condition is satisfied:
\begin{description}
   \item[\ \ A($\delta_h$)\,:] The ideal 
 $\mbox{Ann}_{{\cal D}}\,\delta_h$ of operators annihilating $\delta_h$
is generated by $\mbox{Ann}_{{\cal O}}\,\delta_h$ and operators 
 $Q_1,\ldots, Q_w\in {\cal D}$ of order $1$.
\end{description}
Indeed, because of the relations: 
$Q_i\cdot\delta_hf^{s+1}=(s+1)[Q_i,f]\delta_hf^s$, $1\leq i\leq w$,
we have the following isomorphism:
\begin{eqnarray*}
\frac{{\cal D}[s]\delta_hf^s}{{\cal D}[s](\widetilde{\cal J}_{h,f},f)\delta_hf^s}
&\stackrel{\cong}{\longrightarrow} &(s+1)\frac{{\cal D}[s]\delta_hf^s} {{\cal D}[s]\delta_hf^{s+1}}
\end{eqnarray*}
where $\widetilde{\cal J}_{h,f}\subset {\cal O}$ 
is generated by the commutators 
$[Q_i,f]\in{\cal O}$, $1\leq i\leq w$. Thus $\tilde{b}_f(\delta_h,s)$
may also be defined using the functional equation:
$$b(s)\delta_h f^s\in{\cal D}[s](\widetilde{\cal J}_{h,f},f)\delta_hf^s\ $$
and ${\cal K}={\cal D}[s](\widetilde{\cal J}_{h,f},f)\delta_h f^s/
{\cal D}[s]({\cal J}_{h,f},f)\delta_h f^s$. 
For more details about this condition {\bf A($\delta_h$)}, see \cite{TO3}.

\section{The proofs} \label{proof}

Let us recall that $\tilde{b}_h(s)$ may be defined as the unitary nonzero
polynomial $b(s)\in{\bf C}[s]$ of smallest degree such that:
\begin{equation}
  \label{caracBred}
   b(s)h^s=P(s)\cdot h^{s+1}+\sum_{i=1}^nP_i(s)\cdot\,h'_{x_i}h^s
\end{equation}
where $P(s),P_1(s),\ldots, P_n(s)\in{\cal D}[s]$ (see \cite{12}).

\stl

\begin{rema}{\em 
  The equation (\ref{caracBred}) is equivalent to the following one:
   $$b(s)h^s=\sum_{i=1}^nQ_i(s)\cdot\,h'_{x_i}h^s$$
  where $Q_i(s)\in{\cal D}[s]$ for $1\leq i\leq n$. 
  Indeed, one can prove that 
  $h^{s+1}\in{\cal D}[s](h'_{x_1},\ldots,h'_{x_n})h^s$ 
  {\it i.e.} $h$ belongs to the ideal 
$I={\cal D}[s](h'_{x_1},\ldots, h'_{x_n}) + \mbox{Ann}_{{\cal D}[s]}h^s$.
This requires some computations like in \cite{18} 2.1., 
using that: $h\partial_{x_i}-s h'_{x_i}\in I$, $1\leq i\leq n$.}
\end{rema}

\noindent{\it Proof of Proposition \ref{redBred}}. By semi-continuity
of the Bernstein polynomial, it is enough to prove the assertion for
a reducible germ $f$. Let us write $f=f_1f_2$ where
$f_1,f_2\in{\cal O}$ have no common factor. Assume that $-1$ is not a root of
$\tilde{b}_f(s)$. Then, by fixing $s=-1$ in (\ref{caracBred}), we get:
$$\frac{1}{f}\in\sum_{i=1}^n{\cal D}\frac{f'_{x_i}}{f}+{\cal O}\ \subset \
{\cal O}[1/f_1]+{\cal O}[1/f_2]$$
since $f'_{x_i}/f=f'_{1,x_i}/f_1+f'_{2,x_i}/f_2$, $1\leq i\leq n$. But this
is absurd since $1/f_1f_2$ defines a nonzero element of ${\cal O}[1/f_1f_2]
/{\cal O}[1/f_1] +{\cal O}[1/f_2]$
under our assumption on $f_1,f_2$. Thus $-1$ is a root of $\tilde{b}_f(s)$.
$\square$

\stl

The proofs of the equivalence between 1 and 2 in Theorem \ref{LaiR} and of Theorem \ref{LaiRgen} are based on the following result:

\begin{prop} \label{propmhs}
  Let $f\in{\mathcal O}$ be a nonzero germ such that $f(0)=0$.
   Let  $m$ be a section of a holonomic ${\mathcal D}$-module ${\cal M}$
    without $f$-torsion, and $\ell\in{\bf N}^*$. The following conditions
   are equivalent:
 \begin{enumerate}
   \item \label{ppetite} The smallest integral root of $b_f(m,s)$ is
   strictly greater than $-\ell-1$.
   \item \label{engloc} The ${\mathcal D}$-module $({\mathcal D}m)[1/f]$
    is generated by $mf^{-\ell}$.  
   \item \label{engcoh} The ${\mathcal D}$-module 
   $({\mathcal D}m)[1/f]/{\mathcal D}m$
    is generated by $\stackrel{.}{mf^{-\ell}}$.
   \item \label{eval} The following ${\cal D}$-linear morphism
 is an isomorphism~:
     \begin{eqnarray*}
          \pi_\ell~:~  \frac{{\cal D}[s]mf^s}{(s+\ell){\cal D}[s]mf^s}
             &\longrightarrow & ({\cal D}m)[1/f] \\
             \stackrel{.}{P(s)\cdot\,mf^s} &\mapsto & P(-\ell)\cdot\,mf^{-\ell}
     \end{eqnarray*}
\end{enumerate}
\end{prop}

\begin{demo}
This is a direct generalization of a well known result due to
M. Kashiwara and E. Bj\"ork for $m=1\in {\mathcal O}={\mathcal M}$
(\cite{11} Proposition 6.2, \cite{Bj} Proposition 6.1.18, 6.3.15 \& 6.3.16).

  Let us prove ${\it \ref{ppetite}} \Rightarrow{\it\ref{eval}}$. 
First, we establish that $\pi_\ell$ is surjective. It is enough to see 
that for all $P\in{\cal D}$ and $l\in{\mathbf Z}$~: 
  $(P\cdot\,m)f^l\in{\cal D}mf^{-\ell}$. By using the following relations:
$$\left(\left[\frac{\partial}{\partial x_i}Q\right]\cdot m\right)f^l=\left[\frac{\partial}{\partial
    x_i}f-l\frac{\partial f}{\partial x_i}\right]\cdot\,((Q\cdot m)f^{l-1})$$
  where $1\leq i \leq n$, $Q\in{\cal D}$ and $l\in{\mathbf Z}$, we obtain that
 for all
  $P\in{\cal D}$, $l\in{\mathbf Z}$, there exist $Q\in{\cal D}$ and $k\in{\mathbf Z}$ such that
  $(P\cdot m)f^l=Q\cdot\,mf^k$. Thus, we just have to prove that:
    $mf^k\in{\cal D}mf^{-\ell}$ for $k<-\ell$. 

  Let $R\in{\cal D}[s]$ be a differential operator such that:
  \begin{equation}
    \label{polBer}
     b_f(m,s)mf^s=R\cdot\,mf^{s+1}
  \end{equation}
  and let $k\in{\mathbf Z}$ be such that $k<{-\ell}$. Iterating
  (\ref{polBer}), we get the following identity in ${\cal D}m[1/f,s]f^s$: 
  \begin{equation}
    \label{politr}
   \underbrace{b_f(m,s-{\ell}-k-1)\cdots
   b_f(m,s+1)b_f(m,s)}_{c(s)}mf^s=Q(s)\cdot\,mf^{s-{\ell}-k} 
  \end{equation}
where $Q(s)\in{\cal D}[s]$. By assumption on $\ell$, we have: $c(k)\not=0$. Thus, 
by fixing $s=k$ in (\ref{politr}), we get $mf^k\in{\cal D}mf^{-\ell}$ and 
$\pi_\ell$ is surjective.

\stl

Let us prove the injectivity of $\pi_\ell$. If we fix $P(s)\in{\cal D}[s]$, then
we have the following identity in ${\cal D} m[1/f,s]f^s$:
$$P(s)\cdot\,mf^s=(Q(s)\cdot m)f^{s-l}$$
where $Q(s)\in{\cal D}[s]$ and $l$ is the degree of $P$. Assume that
$\stackrel{.}{P(s)\cdot\,mf^s}\in\ker \pi_\ell$. Thus there exists a non negative
integer $j\in{\mathbf N}$ such that $f^jQ(-\ell)$ annihilates $m\in{\cal M}$. 
In particular: $P(s)\cdot\,mf^s=(s+{\ell})(Q'\cdot m)f^{s-l}$, where $Q'\in{\cal D}[s]$
 is the quotient of the division  of $Q$ by
$(s+{\ell})$. As in the beginning of the proof, we obtain that 
$P(s)\cdot\,mf^s=(s+{\ell})\widetilde{Q}\cdot\,mf^{s-k}$ where $\widetilde{Q}\in{\cal
 D}[s]$
and $k\in{\mathbf N}^*$. From (\ref{polBer}), we get:
$$\underbrace{b_f(m,s-1)\cdots
  b_f(m,s-k+1)b_f(m,s-k)}_{d(s)}P(s)\cdot\,mf^s=(s+{\ell})\widetilde{Q}S\cdot\,mf^s$$
where $S\in{\cal D}[s]$. By division of $d(s)$ by $(s+{\ell})$, we obtain
the identity~:
$$d({-\ell})P(s)\cdot\,mf^s=(s+{\ell})[\widetilde{Q}S+e(s)P(s)]\cdot\,mf^s$$
where $e(s)\in{\mathbf C}[s]$. Remark that $d({-\ell})\not=0$ by assumption on $\ell$. 
Thus $P(s)\cdot\,mf^s\in(s+{\ell}){\cal D}[s]mf^s$, and $\pi_\ell$ is injective. Hence the
condition {\it \ref{ppetite}} implies that
$\pi_\ell$ is an isomorphism.

\stl

Observe that  ${\it \ref{eval}\Rightarrow\ref{engloc}}$ and
  ${\it\ref{engloc}\Leftrightarrow\ref{engcoh}}$ are clear. Thus 
 let us prove 
${\it \ref{engloc}\Rightarrow\ref{ppetite}}$. 
Let $k\in{\mathbf Z}$ denote the smallest integral root of $b_f(m,s)$.
Assume that ${-\ell}>k$. We have the following commutative diagram: 

$$\begin{array}{cccccccc}
         & &      0 & &0 & & &  \\
         & &      \downarrow & & \downarrow &  & &\\
  0 &\rightarrow & {\cal D}[s]mf^{s+1} & \hookrightarrow & {\cal
  D}[s]mf^{s} & \twoheadrightarrow &{{\cal D}[s]mf^{s}}/{{\cal
 D}[s]mf^{s+1}}
  &\rightarrow 0  \\
 & &  \downarrow & & \downarrow & &  \downarrow \upsilon & \\
  0 &\rightarrow & {\cal D}[s]mf^{s+1} & \hookrightarrow & {\cal
  D}[s]mf^{s} & \twoheadrightarrow & {{\cal D}[s]mf^{s}}/{{\cal
         D}[s]mf^{s+1}} & \rightarrow 0 \\ 
  & &          \downarrow u & & \downarrow & & &\\
  & &{\cal D}mf^{k+1} &\stackrel{i}{\hookrightarrow}  & {\cal
         D}mf^k  & &  & \\
   & &        & & \downarrow  & & &\\
     & &       & &0 & &  &
 \end{array}$$

\noindent where $\upsilon$ is the left-multiplication by $(s-k)$. Remark that
 the second column is exact (since 
 ${\it \ref{ppetite}\Rightarrow\ref{eval}}$), that $u$ is surjective,
 and that $i$ is an isomorphism (since $mf^{-\ell}\in{\cal
 D}mf^{k+1}$ generates ${\cal D}m[1/f]$ by assumption).

 After a diagram chasing, one can check that $\upsilon$ is surjective. 
 Thus the ${\cal D}$-module ${\cal
   D}[s]mf^s/{\cal D}[s]mf^{s+1}$ is Artinian, as the stalk of a holonomic 
   ${\cal D}$-Module [indeed, it is the quotient of two sub-holonomic
 ${\cal D}$-Modules which are isomorphic (see \cite{11},
 \cite{K2})]. As a surjective endomorphism of an Artinian module is also injective, 
 $\upsilon$ is injective. But this is absurd since $k$ is a root of $b_f(m,s)$. 
 Hence, ${-\ell}$ is less or equal to the smallest integral root of $f$.
\end{demo}

\begin{rema}{\em
Obviously, the statement does not work for any $\ell\in {\bf Z}$ (take 
$m=1\in{\mathcal O}={\mathcal M}$ and $\ell=-1$). Nevertheless, it
is true for any $\ell\in{\bf C}$ such that for all root $q\in{\bf Q}$
of $b_f(m,s)$, we have: $-\ell-q\not\in {\bf N}^*={\bf N}\cap {\bf C}^*$.}
\end{rema}

 Obviously, Proposition \ref{criterion} is obtained by iterating this result.
 Let us prove Theorem \ref{LaiRgen}.

\stl

\noindent{\it Proof of Theorem \ref{LaiRgen}}. 
If $b'_f(h,s)$ has no strictly negative integral root, 
then ${\cal D}\delta_{h,f}={\cal R}_{h,f}$ (Lemma \ref{relatc}
 and Proposition \ref{propmhs},
using that ${\cal D}\delta_h={\cal R}_h$), and we just 
have to remark that $\delta_{h,f}$ belongs to 
${\cal L}_{h,f}$ when $-1$ is not a
root of $b'_f(h,s)$. Indeed, by fixing $s=-1$ in the
defining equation of $b'_f(h,s)$:
 $$b'_f(h,s)\delta_h f^s \in  {\cal D}[s]({\cal J}_{h,f},f)\delta_hf^s$$
 we get:
$$\delta_h f^{-1}\in
\sum_{1\leq k_1<\cdots <k_p\leq n}^n{\cal D} 
m_{k_1,\ldots,k_p}(h,f)\delta_h f^{-1}+{\cal D}\delta_h\subset{\cal R}_h[1/f].$$ 
Thus $\delta_{h,f}\in
{\cal R}_{h,f}\cong{\cal R}_h[1/f]/{\cal R}_h$ belongs to  
${\cal L}_{h,f}$.

Now let us assume that ${\cal L}_{h,f}={\cal R}_{h,f}$. 
As ${\cal L}_{h,f}\subset {\cal D}\delta_{h,f}$, we also have 
${\cal D}\delta_{h,f}={\cal R}_{h,f}$ {\it i.e.} 
$-1$ is the smallest integral root of  $b_f(\delta_h,s)$  
(Proposition \ref{propmhs}, using the asumption 
${\cal D}\delta_h={\mathcal R}_h$). So let us
prove that $-1$ is not a root of $b'_f(h,s)$, following the
formulation of \cite{Walt} Lemma 1.3. Since 
$\delta_{h,f}\in{\cal L}_{h,f}=\sum_{1\leq k_1<\cdots <k_p\leq n}{\cal D} 
m_{k_1,\ldots,k_p}(h,f)\delta_{h,f}$, 
we have: $1\in {\cal D}{\cal J}_{h,f} + \mbox{Ann}_{\cal D}\,\delta_{h,f}$, 
or equivalently: 
$1\in {\cal D}({\cal J}_{h,f},f)+\mbox{Ann}_{\cal D}\,\delta_h f^{-1}$ 
(using that ${\cal D}f(\delta_h f^{-1})={\cal R}_h$). 
Moreover, as $-1$ is the smallest integral root of $b_f(\delta_h,s)$, 
an operator $P$ belongs to $\mbox{Ann}_{\cal D}\,\delta_h \otimes 1/f$ 
if and only if there exists $Q(s)\in{\cal D}[s]$ 
such that $P-(s+1)Q(s)\in \mbox{Ann}_{{\cal D}[s]}\,\delta_h f^s$
(Proposition \ref{propmhs}). Thus we have:
$${\cal D}[s]={\cal D}[s](s+1,{\cal J}_{h,f},f)+
\mbox{Ann}_{{\cal D}[s]}\,\delta_h f^s.$$ 
In particular, if $(s+1)$ was a factor of $b'_f(h,s)$, we would have:
$$\frac{b'_f(h,s)}{s+1}\in
{\cal D}[s](b'_f(h,s),{\cal J}_{h,f},f)+
\mbox{Ann}_{{\cal D}[s]}\,\delta_h f^s$$
But from the identity (\ref{eqcarbprim}), we have:
$$b'_f(h,s)\in{\cal D}[s]({\cal J}_{h,f},f)+\mbox{Ann}_{{\cal D}[s]}\,\delta_h f^s$$
and this is a defining equation of $b'_f(h,s)$. Thus:
$$\frac{b'_f(h,s)}{s+1}\in {\cal D}[s]({\cal J}_{h,f},f)+
\mbox{Ann}_{{\cal D}[s]}\,\delta_h f^s$$
\noindent In particular, $b'_f(h,s)$ divides
$b'_f(h,s)/(s+1)$, which is absurd. 
Therefore $-1$ is not a root of
$b'_f(h,s)$, and this ends the proof. $\square$

\stl

\noindent{\it Proof of the equivalence between $1$ and $2$ in Theorem \ref{LaiR}}. 
Up to notational changes, the proof is the very same than the previous one.
Assume that $\tilde{b}_h(s)$ has no integral root. On one hand, ${\cal D}\delta_h$ 
coincides with ${\cal R}_h={\cal O}[1/h]/{\cal O}$ by Proposition \ref{propmhs} [take $m=1$ and
${\cal M}={\cal O}$]. On
the other hand, by fixing $s=-1$ in (\ref{caracBred}), we get
$$\frac{1}{h}\in \sum_{i=1}^n {\cal D}\frac{h'_{x_i}}{h} +{\cal O}\subset {\cal O}[1/h]$$
and $\delta_h\in {\cal L}_h$. Hence ${\cal L}_h={\cal R}_h$.

Now let us assume that ${\cal L}_h={\cal R}_h$. As ${\cal L}_h\subset{\cal D}\delta_h \subset {\cal R}_h$, $\delta_h$ generates ${\cal R}_h$. In particular, $-1$ is the only integral root of
$b_h(s)$ by using Proposition \ref{propmhs} (since the roots of $b_h(s)$ are strictly negative).
By the same arguments as in the proof of Theorem \ref{LaiRgen}, one can prove that $-1$ is not
a root of $\tilde{b}_h(s)$. Thus $\tilde{b}_h(s)$ has no integral root, as it was expected $\square$
 
\stl

\begin{rema}{\em
Under the assumption ${\cal D}\delta_{h,f}={\cal R}_{h,f}$, 
we show in the proof of Theorem \ref{LaiRgen} that if  
$\delta_{h,f}$ belongs to ${\cal L}_{h,f}$ then
 $-1$ is not a root of $b'_f(h,s)$. As the reverse
 relation is obvious, a natural question is to know if this 
 assumption is necessary. In terms of reduced Bernstein polynomial,
 does the condition: \textsl{$-1$ is not a root of $\tilde{b}_h(s)$}
 caracterize the membership of $\delta_h$ in ${\cal L}_h$ ?}
\end{rema}

\section{Some remarks}

Let us point out some facts about Theorem \ref{LaiRgen}:

\stl

- The assumption ${\cal D}\delta_h={\cal R}_h$ is necessary. This
appears clearly in the following examples.

\begin{exem}{\em \label{lastexample}
Let $p=1$ and $h=x_1^2+\cdots+x_4^2$. As $b_h(s)=(s+1)(s+2)$,
we have ${\cal D}\delta_h\not={\cal R}_h$ (Proposition \ref{propmhs}).
If $f_1=x_1$, then $b'_{f_1}(h,s)=(s+2)$ by using 
Proposition \ref{propc} 
where as ${\cal L}_{h,f_1}={\cal R}_{h,f_1}$ (Example \ref{exbern}, or 
because ${\cal D}\delta_{f_1}={\cal R}_{f_1}$ and 
$b'_{f_1}(h,s)=(s+3/2)$). 

Now if we take $f_2=x_5$, we have ${\cal L}_{h,f_2}\not={\cal R}_{h,f_2}$
and $b'_{f_2}(h,s)=1$ since:
$$\frac{\stackrel{.}{2}}{x_1^2+\cdots+x_4^2}x_5^s=
\left[\sum_{i=1}^4\frac{\partial}{\partial x_i}x_i\right]\cdot
\frac{\stackrel{.}{1}}{x_1^2+\cdots+x_4^2}x_5^s\ .$$}
\end{exem}
  
 - If $p=1$, this condition ${\cal D}\delta_h={\cal R}_h$ just
 means that the only integral root of $b_h(s)$ is $-1$ 
 (Proposition \ref{propmhs}).

\stl  
  
- The condition ${\cal L}_h={\cal R}_h$ clearly implies 
${\cal D}\delta_h={\cal R}_h$, but it is not necessary; 
see Example \ref{ExLnotR} for instance. An other example
with $p=1$ is given by $h=x_1x_2(x_1+x_2)(x_1+x_2x_3)$
since $b_h(s)=(s+5/4)(s+1/2)(s+3/4)(s+1)^3$.

\stl

- Contrarly to the classical Bernstein polynomial, it may
happen that an integral root of $b_f '(h,s)$ is 
positive or zero (see Example \ref{excalcbprime}, with
$f=x_1$, $h=(x_1^2+\cdots+x_4^2)^\ell$ and $\ell\geq 2$). 
In particular, $1$ is an eigenvalue of the monodromy
acting on $\phi_f{\bf C}_{h^{-1}\{0\}}$. For that reason, we do not have here the analogue of
condition $3$, Theorem \ref{LaiR}. 

\stl

- In \cite{Bud}, the authors introduce a notion of Bernstein 
polynomial for an arbitrary variety $Z$. In the case of hypersufaces, this 
polynomial $b_Z(s)$ coincides with the classical Bernstein-Sato 
polynomial. But it does not seem to us that its integral roots
are linked to the condition 
${\cal L}_{h,f}={\cal R}_{h,f}$. For instance, if 
$h=x_1^2+x_2^2+x_3^2$ and $f=x_4^2+x_5^2+x_6^2$ then
one can check that $b'_f(h,s)=\tilde{b}(f^s,s)=(s+3/2)$;
in particular ${\cal L}_{h,f}={\cal R}_{h,f}$. Meanwhile, 
by using \cite{Bud}, Theorem 5, we get
$b_Z(s)=(s+3)(s+5/2)(s+2)$ if $Z=V(h,f)\subset {\bf C}^6$.

\stl

\end{document}